\documentclass[10pt,leqno]{amsart}

\usepackage[latin1]{inputenc}

\usepackage[a4paper,asymmetric]{geometry}
\usepackage{latexsym, amsmath, amssymb, esint}
\usepackage{bm}
\usepackage{graphics,epsfig,graphicx,graphics}
\usepackage{color, xcolor}
\usepackage{tikz}
\usepackage[hyperpageref]{backref}
\usetikzlibrary{patterns}
\usepackage[english]{babel}
\usepackage[colorlinks=true,urlcolor=blue, citecolor=red,linkcolor=blue,
linktocpage,pdfpagelabels, bookmarksnumbered,bookmarksopen]{hyperref}
\usepackage{hyperref}

\newtheorem{theorem}{Theorem}[section]

\newtheorem{prop}[theorem]{Proposition}

\theoremstyle{definition}

\newcommand{\ep}{\varepsilon}
\newcommand{\RR}{\mathbb{R}}
\newcommand{\rn}{\RR^N}

\newcommand{\de}{\delta}

\newcommand{\pa}{\partial}

\newcommand{\osc}{\rm osc}
\newcommand{\diam}{\mathrm{diam}}
\newcommand{\hd}{\mathrm{hd}}
\newcommand{\dist}{\mathrm{dist}}
\newcommand\Id{{\rm Id}}

\newcommand{\Om}{\Omega}

\numberwithin{equation}{section}

\begin{document}
\title[Almost CMC hypersurfaces]{Quantitative estimates for almost constant mean curvature hypersurfaces}

\author[G. Ciraolo]{Giulio Ciraolo}

\address{G. Ciraolo,  Dipartimento di Matematica ``Federigo Enriques'', Universit\`a degli Studi di Milano, Via Cesare Saldini 50, 20133, Milano - Italy}
\email{giulio.ciraolo@unimi.it}

\maketitle

\begin{abstract}
Alexandrov's soap bubble theorem asserts that spheres are the only connected closed embedded hypersurfaces in the Euclidean space with constant mean curvature. The theorem can be extended to space forms and it holds for more general functions of the principal curvatures.

In this short review, we discuss quantitative stability results regarding Alexandrov's theorem which have been obtained by the author in recent years. In particular, we consider hypersurfaces having mean curvature close to a constant and we quantitatively describe the proximity to a single sphere or to a collection of tangent spheres in terms of the oscillation of the mean curvature. Moreover, we also consider the problem in a non local setting, and we show that the non local effect gives a stronger rigidity to the problem and prevents the appearance of bubbling.

\end{abstract}

%
%

\bigskip




\section{Introduction}
In a series a papers during the sixties \cite{Alex_hyp,Alex2,Alex3}, Alexandrov studied global properties of surfaces with a flavour which is a mix between differential geometry and partial differential equations. One of the most influencing results that he proved is the so called Alexandrov Soap Bubble Theorem, which asserts the following: 

\begin{theorem}[{\sc Alexandrov's Theorem}]
Let $\Omega \subset \mathbb{R}^{n+1}$, $n \geq 1$, be a bounded connected domain with boundary $S=\partial \Omega$ of class $C^2$. Then the mean curvature $H$ of $S$ is constant if and only if $S$ is a sphere.
\end{theorem}
This result was probably expected by the mathematical community, and some attempts and partial results were previously proved \cite{Jellet,Hopf}. Alexandrov's theorem is related to a well-known conjecture of H. Hopf \cite{Hopf}:

\begin{theorem}[{\sc Hopf Conjecture}] 
Let $S$ be an immersion of an oriented, closed hypersurface with constant mean curvature $H \neq 0$ in $\mathbb{R}^{n+1}$.
Must $S$ be the standard embedded $n$-sphere?
\end{theorem}

Since $S=\partial \Omega$ is the boundary of a bounded open set, then the hypersurface $S$ considered in Alexandrov's theorem is embedded and hence Alexandrov proved that Hopf's Conjecture is true under the assumption of embedness. This assumption is optimal as it is was proved by Wente in \cite{We}, where he gave a counterexample of an immersed hypersurface with constant mean curvature which is not a sphere in $\mathbb{R}^3$, nowadays called the Wente's tori (see also \cite{HTY} for a counterexample in $\mathbb{R}^4$). 

Removing the embedness assumption is possible by assuming other assumptions on the hypersurface. For instance, it is a result of Hopf \cite{Hopf} that an immersed and simply connected hypersurface in $\mathbb{R}^3$ with constant mean curvature is a sphere. Other results have been obtained by Barbosa and Do Carmo \cite{BarbosaDoCarmo} under the assumption that the hypersurface is stable.

In this paper we consider an embedded hypersurface $S=\partial \Omega$ of $\mathbb{R}^{n+1}$ having the mean curvature close to a constant and we are interested in quantifying the proximity of $S$ to some \emph{special} configuration. 

At a first glance, one could think that it is reasonable to have that if $H_S$ is close to a constant, then $S$ is close to a sphere. This is not true, since \emph{bubbling} may appear: it is possible to construct counterexamples showing that an almost constant mean curvature hypersurface is close to an array of tangent spheres in any $C^k$ norm (see for instance \cite{Butscher,Butscher_Mazzeo} and \cite{K1,K2}). 

More precisely, the examples available in literature show that $S$ is a small deformation of almost tangent spheres which are connected by very small necks. In these necks the almost umbilicality of the surfaces is completely lost, since the mean curvature turns to be close to a constant even if the absolute value of the principal curvatures becomes arbitrarily large.
This phenomenon suggests that in order to have the proximity to a single sphere, one has to introduce some condition that prevents the blow-up of the principal curvatures. 

Quantitatively describing the bubbling and/or the proximity to a single sphere is the main goal of this note, in which we review some recent results obtained by the author and his collaborators as well as other related results.

\section{Quantifying the bubbling} \label{section_bubbling}
In this section we give a quantitative description of bubbling for almost constant mean curvatures. These results were obtained by the author and Maggi in \cite{CM}.

Let $\Omega \subset \mathbb{R}^{n+1}$ be a bounded connected open set, with boundary $\partial \Omega$ of class $C^2$. We denote by $H= k_1 + \ldots + k_n$ the mean curvature of $\partial \Omega$ (not normalized by $n$), and by $P(\Omega)$ and $|\Omega|$ the perimeter and the $(n+1)$-Lebesgue measure of $\Omega$, respectively.

We begin with some simple consideration. We first notice that if $\partial \Omega$ has constant mean curvature then $H=H_0$, where 
\begin{equation}\label{H0_def}
H_0 = \frac{n P(\Omega)}{(n+1) |\Omega|}\,.
\end{equation}
Indeed, \eqref{H0_def} follows from the tangential divergence theorem and the divergence theorem under the assumption that $H$ is constant:
$$
 n P(\Omega) =\int_{\partial \Omega} H_0 x \cdot \nu d\mathcal{H}^n = H_0 \int_{\partial \Omega} x \cdot \nu =  (n+1) |\Omega| H_0 \,.
$$
We use this observation to introduce a scale invariant deficit, which will measure the distance of $H$ from the constant $H_0$:
\begin{equation} \label{deficit_delta}
\delta (\Omega ) = \Big{\|} \frac{H}{H_0} - 1\Big{\|}_{C^0(\partial \Omega)} \,.
\end{equation}
This is a scale invariant quantity, in the sense that $\delta (\lambda \Omega)=\delta (\Omega)$ for any $\lambda >0$ and $\delta(\Omega)=0$ if and only if $H$ is constant (i.e. $H=H_0$). Clearly, by Alexandrov Theorem, $\delta (\Omega) = 0$ if and only if $\Omega$ is a ball.

Since the deficit $\delta(\Omega)$ is scale invariant, we can assume that $H_0=n$ (hence in case $\delta(\Omega)=0$ then $\Omega=B_1$). In the following theorem, we just assume an upper bound on the perimeter, and we quantitatively describe what happens when $\delta(\Omega)$ is close to zero. Since we only assume a bound on the perimeter, we have to deal with possible bubbling and the bound on the perimeter is needed in order to bound the number of possible bubbles. More precisely, we introduce a scale-invariant quantity 
$$
\mathcal{Q}(\Omega) = \left( \frac{P(\Omega)}{P(B)} \right)^{n+1} \left(\frac{|B_1|}{|\Omega|} \right)^n,
$$
and notice that, thanks to the isoperimetric inequality, one always have $\mathcal{Q}(\Omega) \geq 1$ and that
$$
\mathcal{Q}( \textmd{\rm a union of } L \  \textmd{\rm disjoint balls of equal radii} ) = L
$$
$\forall L \in \mathbb{N}$, $L\geq 1$. As we will show, the functional $\mathcal{Q}$ counts how many disjoint balls of radius $n/H_0$ will approximate $\Omega$. Indeed we have the following theorem which was proved in 

\begin{theorem}\label{thm_CM}
  Given $n,L\in\mathbb N$ with $n\geq 2$ and $L\geq 1$, and $a\in(0,1]$, there exists a positive constant $c(n,L,a)>0$ with the following property. If $\Omega$ is a bounded connected open set with $C^2$-boundary in $\mathbb R^{n+1}$ such that $H>0$ and
  \[
  H_0=n\,,\qquad P(\Omega)\leq (L+1-a)\,P(B)\,,\qquad \de(\Omega)\leq c(n,L,a)\,,
  \]
  then there exists a finite family $\{B_{z_j,1}\}_{j\in J}$ of mutually disjoint balls with $\#\,J\leq L$ such that if we set
  \[
  G=\bigcup_{j\in J}B_{z_j,1}\,,
  \]
  then
  \begin{eqnarray}
  \label{main thm Omega meno G}
  \frac{|\Omega\Delta G|}{|\Omega|}&\leq&  C(n)\,L^2\,\de(\Omega)^\alpha\,,
  \\
  \label{main thm perimetri stima}
    \frac{|P(\Omega)-\#\,J\,P(B)|}{P(\Omega)}\,&\leq& C(n)\,L^2\,\de(\Omega)^\alpha\,,
  \\
  \label{main thm onesided hd}
  \frac{\max_{x\in \partial G}\dist(x,\partial\Omega)}{\diam(\Omega)}&\leq& C(n)\,L\,\de(\Omega)^\alpha\,,
  \\\label{main thm hd stima}
    \frac{\hd(\partial\Omega,\partial G)}{\diam(\Omega)}&\leq& C(n)\,L^{3/n}\,\de(\Omega)^{\alpha/4n^2(n+1)}\,.
  \end{eqnarray}
  Moreover, there exists an open subset $\Sigma$ of $\partial G$ and a function $\psi:\Sigma\to\mathbb R$ with the following properties. The set $\pa G\setminus\Sigma$ consists of at most $C(n)\,L$-many spherical caps whose diameters are bounded by $C(n)\,\de(\Omega)^{\alpha/4(n+1)}$. The function $\psi$ is such that $(\Id+\psi\,\nu_G)(\Sigma)\subset\partial\Omega$ and
  \begin{gather}\label{stima psi intro}
  \|\psi\|_{C^{1,\gamma}(\Sigma)}\le C(n,\gamma)\,,\qquad\forall\gamma\in(0,1)\,,
  \\
  \label{stima psi intro 2}
  \frac{\|\psi\|_{C^0(\Sigma)}}{\diam(\Omega)}\leq C(n)\,L\,\de(\Omega)^\alpha\,,\qquad \|\nabla\psi\|_{C^0(\Sigma)}\leq C(n)\,L^{2/n}\,\de(\Omega)^{\alpha/8n(n+1)}\,,
  \\\label{stima bordo omega meno immagine di sigma}
  \frac{\mathcal H^n(\partial\Omega\setminus(\Id+\psi\,\nu_G)(\Sigma))}{P(\Omega)}\leq C(n)\,L^{4/n}\,\delta(\Omega)^{\alpha/4n(n+1)}\,,
  \end{gather}
  where $(\Id+\psi\,\nu_G)(x)=x+\psi(x)\,\nu_G(x)$ and $\nu_G$ is the outer unit normal to $G$.
\end{theorem}
From a qualitative point of view,  Theorem \ref{thm_CM} has the consequence that examples available in literature on almost constant mean curvature are actually the only possible ones which are not close to a single sphere, and this qualitative information is optimal. 

The quantitative estimates in Theorem \ref{thm_CM} are clearly the main result of the theorem. Indeed, not only they describe quantitatively the appearing of bubbling (although arguably in a non-sharp way), but they can also be used to describe capillarity droplets, as we will describe later.
Moreover, these estimates have a simple and interesting consequence, which is described in the following proposition (see \cite{CM}[Proposition 1.1]).

\begin{prop}\label{proposition referee 1}
  Under the same assumptions of Theorem \ref{thm_CM}:

  \begin{enumerate}
    \item[(i)] if $\#\,J\geq 2$, then for each $j\in J$ there exists $\ell\in J$, $\ell\neq j$, such that
  \begin{equation}
      \label{sjsjprimo tangenti}
    \frac{\dist(\pa B_{z_j,1},\pa B_{z_\ell,1})}{\diam(\Omega)}\leq C(n)\,\de(\Omega)^{\alpha/4(n+1)}\,,
  \end{equation}
  that is to say, each ball in $\{B_{z_j,1}\}_{j\in J}$ is close to be tangent to another ball from the family;
\\
  \item[(ii)] if there exists $\kappa \in(0,1)$ such that
  \begin{equation}
    \label{stima densita fuori intro}
    |B_{x,r}\setminus\Omega|\geq \kappa \,|B|\,r^{n+1}\,,\qquad\forall x\in\partial\Omega\,,r<\kappa\,,
  \end{equation}
  and $\de(\Omega)\leq c(n,L,\kappa)$, then $\#\,J=1$, that is, $\Omega$ is close to a single ball.
  \end{enumerate}
\end{prop}

Item (ii) in Proposition \ref{proposition referee 1} gives a smallness criterion for proximity to a single ball which is very weak: indeed only an exterior uniform volume estimate is needed in order to avoid the appearing of bubbling. As we will show later, this criterion is also important because it can also be applied to describe the shape of local minimizers in capillarity problems.

Now we describe the strategy of the proof of Theorem \ref{thm_CM}. To clarify the exposition, we start by giving the proof of Alexandrov Theorem given by A. Ros in \cite{Ros}, which is the starting point for our qualitative analysis.

The proof in \cite{Ros} is based on the \emph{Heintze-Karcher} inequality, which asserts that if $H>0$ then
\begin{equation} \label{ineq_HK}
\int_{\partial \Omega} \frac{n}{H} d\mathcal{H}^n \geq (n+1) |\Omega| \,.
\end{equation}
This inequality can be proved by introducing an auxiliary problem, the torsion problem for $\Omega$, i.e. considering 
$$
\begin{cases}
\Delta f = 1 & \text{ in } \Omega \,, \\
f= 0 & \text{ on } \partial \Omega \,,
\end{cases}
$$
and by showing that, by using Reilly's inequality (see \cite{Reilly}), the equality in \eqref{ineq_HK} is attained if and only if 
$$
\nabla^2 f = \frac{1}{n+1} \Id \quad \text{ in }  \Omega \,,
$$
and 
$$
|\nabla f| = \frac{n}{(n+1) H_0} \quad \text{ on } \partial \Omega  \,.
$$
The proof of Theorem \ref{thm_CM} can be summarized as follows.

\emph{Step 1 - Uniform bounds for the torsion potential.} We use classical symmetrization results and a $P$-function approach to prove that
  \begin{eqnarray}
    \label{f C0}
    \|f\|_{C^0(\Om)}&\le & C\,|\Om|^{2/(n+1)}\,,
    \\
    \label{f Lip}
    \|\nabla f\|_{C^0(\Om)}&\le& C_0\,|\Om|^{1/(n+1)}\,,
    \\
    \label{f L2 estimates}
    \|\nabla^2f\|_{L^2(\Om)}&\le&|\Om|^{1/2}\,.
  \end{eqnarray}

\emph{Step 2 - Quantitative analysis of Ros' proof.} By looking at Ros' proof by a quantitative viewpoint and making use of \emph{Step 1}, we obtain the following two inequalities:
\begin{eqnarray}
\label{quant1}
C(n)\,|\Omega|\,\de(\Omega)^{1/2}&\geq&\int_\Omega\Big|\nabla^2 f-\frac{\Id}{n+1}\Big|\,,
\\
\label{quant2}
C(n)\,\Big(\frac{n}{H_0}\Big)^2\,P(\Omega)\,\de(\Omega)&\geq&\int_{\pa\Omega}\Big|\frac{n/H_0}{n+1}-|\nabla f|\Big|^2\,,
\end{eqnarray}
where the second estimate holds if $\de(\Omega)\leq1/2$, and where $\nabla f=|\nabla f|\,\nu_\Omega\neq 0$ on $\pa\Omega$.

\emph{Step 3 - Localize the estimates in Step 2.} We localize the estimates \eqref{quant1} and \eqref{quant2} by considering convolutions $f_\ep $ of the torsion potential $f$. This implies that, on the smaller set $\Omega_\ep = \{x \in \Omega :\ \dist(x,\partial \Omega) > \ep\}$, we have the following pointwise bounds:
  \begin{eqnarray}\label{feps nabla limitato}
  \|\nabla f_\ep\|_{C^0(\Om_\ep)}&\le& C_0(n)\,,
  \\
  \label{feps f in C0}
  \|f_\ep-f\|_{C^0(\Om_\ep)}&\le& C_0(n)\,\ep\,,
  \\
  \label{D2feps identita C0}
  \big\|\nabla^2f_\ep-\frac{\Id}{n+1}\big\|_{C^0(\Om_\ep)}&\le& C(n)\,\eta(\Om)^\alpha\,,
  \end{eqnarray}
  \vspace{-0.5cm}
  \begin{eqnarray}
  \label{D2feps positivo}
  \|\nabla^2f_\ep\|_{C^0(\Om_\ep)}\le C(n)\,,\qquad \nabla^2f_\ep(x)\ge\frac{\Id}{2(n+1)}\,,\qquad\forall x\in\Om_\ep\,,
  \end{eqnarray}
 where $\eta(\Omega)$ is called the \emph{Heintze-Karcher} deficit and it is defined by 
$$
\eta(\Omega) = 1 - \frac{(n+1)|\Omega|}{\int_{\partial \Omega} \frac{n}{H} } \,. 
$$

\emph{Step 4 - Prove the existence of approximating balls}. A crucial step is to prove that there exist $\{B_{x_j,s_j}\}_{j\in J}\subset\{B_{x_i,r_1^i}\}_{i\in I}$ such that 
\begin{equation}
     \label{main thm sj-1}
  \frac{\max_{j\in J}\,|s_j-1|}{\diam(\Om)}\le C(n)\,|\Om|\,\de(\Om)^\alpha\,,
\end{equation}
and, if $G^*=\bigcup_{j\in J}B_{x_j,s_j}$ (so that $G^*\subset\Om$ by construction), then
\begin{gather}
  \label{main thm Omega meno G^*}
  \frac{|\Om\setminus G^*|}{|\Om|}\le  C_1(n)\,\diam(\Om)\,|\Om|\,\de(\Om)^\alpha\,,
  \\
  \label{main thm perimetri stima proof}
    \frac{|P(\Om)-\#\,J\,P(B)|}{P(\Om)}\,\le C(n)\,\diam(\Om)\,|\Om|\,\de(\Om)^\alpha\,,
    \\
    \label{starstar}
    \#\,J\le L\,,\qquad \#\,J\le C(n)\,|\Om|\,.
\end{gather}
This step is proved by exploiting the estimates in Step 3 and obtaining information on the set $A_\ep = \{f_|ep <- 3 \rho\}$, where $\rho= C_0(n) \ep$. In particular, we can show that for $\ep$ small enough we have
$$
\{f < - 4\rho \} \subset A^\ep \subset \{f < -2 \rho\}
$$
and that each connected component of $A^\ep$  is convex as well as other quantitative information on the $A$. Then, by carefully exploiting Pohozaev's identity in a quantitative way, we obtain the estimates \eqref{main thm sj-1}--\eqref{starstar}.
 
\emph{Step 5 - Quantitative proximity of $G^*$ to $\Omega$ and quantitative estimates for $G$}. 
By refining the estimates on $\nabla f_\ep$, in particular by proving that
$$
\Big|\nabla f_\ep(x_0)-\frac{(x_0-x_i)}{n+1}\Big|\le C\,\de(\Om)^\alpha\,,
$$
we can show that
\begin{eqnarray}
  \label{main thm onesided hd G*}
  \frac{\max_{x\in \pa G^*}\dist(x,\pa\Om)}{\diam(\Om)}\,\le C(n)\,\de(\Om)^\alpha\,.
\end{eqnarray}
Finally, starting from $G^*$, we are able to construct the set $G$ as in the statement of the theorem which has the desired properties.

\emph{Step 6 - Almost all of $\partial \Omega$ is a normal graph on $\partial G$.} 
In order to prove this step, the main idea is to use the area excess regularity criterion of Allard. In particular, at points $x\in\pa\Om$ which are sufficiently close to $\pa G$, we need to quantify the size of $\mathcal{H}^n(\partial\Om\cap B_{x,r})$, where $r$ is \emph{small} and proportional to a suitable power of $\delta(\Om)$. This is done by carefully partitioning $\mathbb R^{n+1}$ into suitable polyhedral regions associated to the balls $B_{z_j,1}$, and by then performing inside each of these regions a calibration type argument with respect to the corresponding ball $B_{z_j,1}$ (see \cite{CM}[Proof of Theorem 1.1 - Step six]).

Thanks to this argument, and many others, we are able to parameterize a large portion of $\pa\Om$ over a large portion of $\pa G$ and obtain the estimates \eqref{stima psi intro}, \eqref{stima psi intro 2} and \eqref{stima bordo omega meno immagine di sigma}, which completes the proof of Theorem \ref{thm_CM}.


\section{Proximity to a single ball}
As we have seen in Section \ref{section_bubbling},  an almost constant mean curvature hypersurface is not necessarily close to a sphere and instead bubbling may appear. In this section we describe some quantitative results where the hypersurface is close to a sphere. In order to do that, one has to introduce some further assumption on the hypersurface which prevents bubbling, as we are going to describe below.

Quantitative results for almost constant mean curvature hypersurface have been largely studied under the assumption that the domain is convex: if $\Omega$ is an ovaloid, the problem was studied by Koutroufiotis \cite{Kou},  Lang \cite{La} and Moore \cite{Moo}. Other stability results can be found in Schneider \cite{Schn} and Arnold \cite{Ar}. These results were improved by Kohlmann in \cite{Kol} where he proved an explicit H\"{o}lder type stability in \eqref{stab_est} below. It is clear from Theorem \ref{thm_CM} that the assumption that $\Omega$ is convex forces the domain to be close to a ball as the mean curvature goes to a constant. Moreover, from the quantitative point of view, the estimates obtained in the papers cited above are not sharp.

These estimates where largely improved by the author and Vezzoni in \cite{CV_JEMS} under the weaker assumption that the hypersurface satisfies a touching ball condition of fixed radius, i.e. by assuming that there exists $\rho>0$ such that for any $x \in \partial \Omega$ there exists a ball of radius $\rho$ touching $\partial \Omega$ at $x$ from inside $\Omega$ and a ball of radius $\rho$ touching $\partial \Omega$ at $x$ from outside $\Omega$. The main result in \cite{CV_JEMS} is a quantitative stability estimate for Alexandrov's theorem in terms of the oscillation of the mean curvature
$$
\mathrm{osc}(H)=\max_{p\in S}H(p)-\min_{p\in S}H(p)\, .
$$

\begin{theorem} [\cite{CV_JEMS}]\label{thm_CV}
	Let $S= \partial \Omega$ be an $n$-dimensional, $C^2$-regular, connected, closed hypersurface embedded in $\mathbb{R}^{n+1}$, with $\Omega \subset \mathbb{R}^{n+1}$ a bounded domain satisfying a touching ball condition of radius $\rho$.
	There exist constants $\varepsilon, C>0$ such that if
	\begin{equation}
	\mathrm{osc}(H) \leq\varepsilon\, ,
	\end{equation}
	then there are two concentric balls $B_{r}$ and $B_{R}$ such that
	\begin{equation}
	B_r\subseteq \Omega\subseteq B_R\, ,
	\end{equation}
	and 
	\begin{equation}\label{stab_est}
	R-r\leq C\mathrm{osc}(H) \, .
	\end{equation}
	The constants $\varepsilon$ and $C$ only depend on $n$ and an upper bound on $\rho^{-1}$ and on $|S|$.
	
	Moreover, $S$  is diffeomorphic to a sphere and there exists a $C^1$ map 
	$$
	F=Id + \Psi \nu: \partial B_{r_i} \to S
	$$ 
	s.t.
\begin{equation} \label{stima_C1_JEMS}
\|\Psi\|_{C^1(\partial B_{r_i})} \leq C \ {\rm osc} \, H \,,
\end{equation}
where $C$ depends only on $n$ and an upper bound on $\rho^{-1}$ and on $|S|$.
\end{theorem}

In view of the above in Section \ref{section_bubbling}, the assumption on the touching in Theorem \ref{thm_CV} is the condition which prevents bubbling, since it gives a constraint on the principal curvatures (if bubbling appears then one has curvatures going to infinity in the necks connecting the tangent almost spheres). 

The proof of Theorem \ref{thm_CV} is based on a quantitative analysis of the original proof of Alexandrov. 
Alexandrov's theorem was important not only for the result itself, but also for the technique that Alexandrov used to prove it: the method of moving planes (MMP). This method, as well as Alexandrov's theorem, still have a great influence nowadays on the research in geometric analysis and partial differential equations (see the reviews \cite{Brezis} and \cite{CirRoncBP}). 

The proof of Theorem \ref{thm_CV} is based on a quantitative analysis of the method of moving planes. This type of analysis was first done in \cite{ABR} (see also \cite{CMS, CMS2, CirMagnVespri}). As a first hint, one has to replace qualitative tools like maximum principles by quantitative ones as Harnack's inequality. In order to properly describe a sketch of the proof, we recall how the method of moving planes given in the proof by Alexandrov works.

For any fixed direction $\omega$ we consider the hyperplanes $\pi_\lambda = \{x \cdot \omega = \lambda \}$ orthogonal to $\omega$. Since $\Omega$ is bounded, we have that $\Omega \cap \pi_\lambda= \emptyset$ for $\lambda$ large. We decrease the value of $\lambda$ until we find a hyperplane $\pi_{\lambda_0} $ which is tangent to $\partial \Omega$. Since $\Omega$ is bounded and smooth, for $\lambda< \lambda_0$ and $\lambda$ close to $\lambda_0$ we have that the reflection of 
$$
\Omega_\lambda = \{x \in \Omega:\ x \cdot \omega > \lambda\}
$$
is contained in $\Omega$. We denote by $\mathcal{R}_\lambda (\Omega_\lambda)$ the reflection of $\Omega_\lambda$ about $\pi_\lambda$. By continuing decreasing $\lambda$, we arrive at a critical position given by
$$
\lambda_* = \inf \{\lambda \in \mathbb{R} :  \ \mathcal{R}_\mu (\Omega_\mu) \subset \Omega \textmd{ for any } \mu \in (\lambda, \lambda_0) \} \,,
$$
where the reflection of $\Omega_{\lambda_*}$ is tangent to $\Omega$. This critical position may occur in two different ways: 

(i) $\mathcal{R}_{\lambda_*} (\Omega_{\lambda_*})$ is tangent to $\Omega$ at a point $p\in \partial \Omega \setminus \pi_{\lambda_*}$;

(ii)   $\mathcal{R}_{\lambda_*} (\Omega_{\lambda_*})$ is tangent to $\Omega$ at a point $q\in \partial \Omega \cap \pi_{\lambda_*}$.

Since $p$ and $q$ are tangency points, the tangent spaces of $\partial \Omega$ and $\mathcal{R}_{\lambda_*} (\Omega_{\lambda_*})$ at $p$ and $q$ coincide. Hence,
In both cases, we can locally write $\partial \Omega$ and $\mathcal{R}_{\lambda_*} (\Omega_{\lambda_*})$ as graphs of function on the tangent space at $p$ or $q$, and we find two functions $u^1,u^2: E \to \mathbb{R}$ which parametrize $\partial \Omega$ and $\mathcal{R}_{\lambda_*} (\Omega_{\lambda_*})$, respectively. The set $E \subset \mathbb{R}^n$ is a subset of the tangent space and it is a ball in case (i) and half a ball in case (ii).

Hence $u^1$ and $u^2$ satisfy the mean curvature equation
$$
 \mathcal{L}(u^i):= {\rm div} \left(\frac{\nabla u^i}{\sqrt{1+|\nabla u^i|^2}} \right) = n H \,,
$$
$i=1,2$ and by construction we have that 
$$
u^1-u^2 \geq 0 \,.
$$
Since $H$ is constant and $\nabla u^i$ are bounded in $E$ (up to choosing a smaller set), then the difference $u^1-u^2$ is nonnegative and it satisfies a linear elliptic equation in $E$
\begin{equation} \label{Lu1u2}
L(u^1-u^2)=\mathcal{L}(u^1) -\mathcal{L}(u^2) = 0 \,.
\end{equation}
In case (i), we have that $E=B_r$ and $u^1(O) -u^2(O) =0$ and by the strong maximum principle we conclude that $u^1 - u^2 = 0$ in $E$. In case (ii) we have that $E=B_r^+$ (half ball) and $u^1(O)-u^2(O) =0$, where now $O \in \partial E$; now the conclusion follows from Hopf's boundary point lemma and we conclude that, again,  $u^1 - u^2 = 0$ in $E$.

Hence, we have proved that the set of tangency points is both closed and open, and then $\partial \Omega \cap \{x \in \Omega:\ x \dot \omega < \lambda_*\}$ and $\mathcal{R}_{\lambda_*} (\Omega_{\lambda_*})$ must coincide, i.e. $\Omega$ is symmetric with respect to the direction $\omega$. 

Since the direction $\omega$ is arbitrary then $\Omega$ is symmetric with respect to any direction and then it is easy to show that it has a radial symmetry. Moreover, by construction of the method of moving planes we also obtain that $\Omega$ is simply connected and the condition that $H$ is constant forces $\Omega$ to be a ball. This completes the proof of Alexandrov's theorem.

Now we briefly review this proof from a quantitative point of view. When we apply the method of moving planes, we still arrive at the two possible critical positions (i) and (ii), we still find two ordered solutions $u^1 \geq u^2$, but since $H$ is not constant \eqref{Lu1u2} is replaced by 
\begin{equation} \label{Lu1u2_quant}
L(u^1-u^2)= f
\end{equation}
where $\|f\|_{C^0} \leq \osc H$. Hence, in case (i) we can apply Harnack's inequality and obtain that there exists $C$ such that 
\begin{equation} \label{harnack1}
\|u^1-u^2\|_{C^0(B_{r/2}(O))} \leq C \osc H \,,
\end{equation}
and by interior elliptic regularity estimates we obtain
\begin{equation} \label{harnack2}
\|u^1-u^2\|_{C^1(B_{r/4}(O))} \leq C \osc H \,.
\end{equation}
Here, the touching ball condition plays an important role. Indeed, thank to this condition, the size of $r$ can be estimated in terms of $\rho$, for instance $r=\rho/2$ is fine. This condition gives also an upper bound on the gradients of $u^1$ and $u^2$, which is crucial in order to apply elliptic estimates. In this review, we only describe the quantitative estimates in case (i); case (ii) is analogous but needs some more technicality. For instance, \eqref{harnack1} must be replaced by a quantitative version of Hopf's boundary point lemma.

We notice that \eqref{harnack2} not only says that the graphs of $u^1$ and $u^2$ are close, but also that the normals to $\partial \Omega$ and $\mathcal{R}_{\lambda_*} (\Omega_{\lambda_*})$ are close in a neighborhood of the tangent point $p$. These two informations are crucial in order to propagate the smallness information in \eqref{harnack2} all around the connected component $\Sigma$ of the reflected cap $\mathcal{R}_{\lambda_*} (\Omega_{\lambda_*})$ which contains the tangency point $p$. This can be achieved by a chain of Harnack's inequality performed with balls of fixed size (thanks to the touching ball condition) all around the reflected cap. This argument is true up to a small error (bounded in terms of $\osc H$) due to the fact that the tangent hyperplanes to $\partial \Omega$ and to $\Sigma$ may not coincide (they coincide only at the first step when the two hypersurfaces are tanget), and we have to take care of this difficulty in order to locally write $\partial \Omega$ and $\Sigma$ as graphs of function on the tangent space of a suitable point in $\Sigma$.

We also emphasize that a Harnack's chain of inequality may be applied up to a fixed distance from the hyperplane $\pi_{\lambda_*}$. When we are close to $\pi_{\lambda_*}$ we have to use quantitative Carleman's type estimates. 

After this argument, we are able to show that the domain $\Omega$ is \emph{almost} symmetric with respect to $\pi_{\lambda_*}$. Here, the word \emph{almost} means the following:

\begin{theorem}[Theorem 4.1 in \cite{CV_JEMS}] \label{thm_approx_1direc}
There exists a positive constant $\ep$ such that if
$$
{\rm osc} {\rm H} \leq \ep,
$$
then for any $p \in\Sigma$ there exists $\hat p\in \hat\Sigma$ such that
\begin{equation}\label{bound on dist}
|p-\hat p| + |\nu_p-\nu_{\hat p}| \leq C\, \osc H  . 
\end{equation}
Here, the constants $\ep$ and $C$ depend only on $n$, $\rho$, $|S|$ and do not depend on the direction $\omega$.
\end{theorem}

Once we have the approximate symmetry in any direction, we obtain the approximate radial symmetry by using the following argument. 

We apply the method of moving planes in the direction of the coordinate axes $e_1,\ldots,e_{n+1}$ and we find $n+1$ critical position and critical hyperplanes $\pi_1,\ldots,\pi_{n+1}$.  Let 
$$
\mathcal{O} = \bigcap_{i=1}^{n+1} \pi_i \,.
$$
The point $\mathcal O$ will serve as approximate center of symmetry. Thanks to the approximate symmetry in the directions $e_1,\ldots,e_{n+1}$, the reflection $\mathcal{R}$ about $\mathcal O$ is given by 
$$
\mathcal R = \mathcal R_{\pi_{n+1}} \circ \cdots \circ \mathcal  \mathcal R_{\pi_{1}} \,,
$$
where $\mathcal R_{\pi_{i}}$ denotes the reflection about the critical hyperplane $\pi_i$ in the direction $e_i$, $i=1,\ldots, n+1$, and the approximate symmetry about $\mathcal O$ can be estimated by using Theorem \ref{thm_approx_1direc} (iterated $(n+1)$-times).

By using again the method of moving planes, we can show that every critical hyperplane in a generic direction $\omega$ is close (in a quantitative way) to $\mathcal O$, and then we can prove the quantitative estimate on $r_e-r_1$ in Theorem \ref{thm_CV}.

The rest of the proof of Theorem \ref{thm_CV} (in particular the $C^1$ estimate on $\Psi$) can proved in two steps: we first find a Lipschitz bound on $\Psi$ by using the method of moving planes, and then we can refine the estimate by using elliptic regularity. This completes the proof of Theorem \ref{thm_CV}.

\medskip

We mention that there are other quantitative results estimating the proximity of almost constant mean curvature hypersurfaces to a single sphere. 

In \cite{KrummelMaggi} Krummel and Maggi prove a sharp stability estimate for almost constant mean curvature hypersurfaces by proving a quantitative version of Almgren's isoperimetric principle. In this context, Almgren's isoperimetric principle asserts the following: if $\Omega$ is a bounded open set with smooth boundary in $\mathbb{R}^{n+1}$, then $H_{\partial \Omega} \leq 1$ implies that $P(\Omega) \geq P(B_1)$, where the equality case is attained if and only if $\Omega$ is a ball of radius one.
Thanks to a sharp quantitative version of this principle, the authors are able to prove sharp stability estimates of proximity to a single ball just by assuming that the perimeter of $\Omega$ is strictly less than two times the perimeter of a ball. Hence, in this case, $P(\Omega) < 2 P(B_1)$ is the assumption that prevents bubbling.

Another approach to quantitative estimates of proximity to a single ball, which is based on integral identities, can be found in a series of papers by Magnanini and Poggesi \cite{magnanini_poggesi1,magnanini_poggesi2,magnanini_poggesi3}. Here, the assumptions that prevents bubbling are bounds on the interior and exterior touching ball condition and on the diameter of $\Omega$. Moreover, the stability result is given in terms of a different deficit.

\section{Quantitative results in other settings}
\subsection{Bubbling for crystals} There are many situations of physical and geometric interest where the Euclidean norm is replaced by another norm in $\mathbb{R}^n$ (see for instance \cite{Tay78}). In particular, one can define the \emph{anisotropic perimeter} by considering an anisotropic surface energy of the form
\begin{equation}\label{anis_surf_energy}
P_F(\Om)=\int_{\pa \Om}F(\nu) d\mathcal{H}^n (x) \,,
\end{equation}
where $F$ is a norm in $\mathbb R^n$, i.e.  a one-homogeneous convex function in $\mathbb R^n$.
A set $E$ is \emph{Wulff shape} of $H$ if there exist $t>0$ and $x_0\in\rn$ such that
$$
E=\{x\in\rn\ :\ F_0(x-x_0)\le t\},
$$
where $F_0$ denotes the dual norm of $F$. It is clear that if $F$ is the Euclidean norm then $P_F(\Om)$ is the usual perimeter of $\Omega$.

Starting from \eqref{anis_surf_energy}, the study of critical points of $P_F(\cdot)$ for volume-preserving variations leads to the condition $H_F(\Omega) = C$, where $H_F(\Omega) $ denotes the anisotropic mean curvature of $\Omega$. 

\begin{theorem}[Anisotropic Alexandrov's Theorem]\label{AleksandrovThm}
Let $F$ be a norm of $\mathbb R^n$ of class $C^2(\rn \setminus \{O\})$ such that $H^2$ is uniformly convex, and let $\partial \Omega$ be a compact hypersurface without boundary embedded in Euclidean space of class $C^2$.
If $H_F(\Omega) $ is constant for every $x\in \pa \Omega$ then $\Omega$ has the Wulff shape of $H$.
\end{theorem}
\noindent A proof of the previous result can be found in \cite{BianchiniCiraoloSalani}[Appendix B], \cite{16BCS} and \cite{17BCS}.

A relevant quantitative study of the anisotropic Alexandrov theorem was done in \cite{DMMN}. In that paper, the authors generalize the results of \cite{CM} to the anisotropic setting and also considered an $L^2$-deficit.

\medskip

\subsection{Quantitative estimates in space forms}

The proximity to a single ball in \cite{CV_JEMS} has been extended by the author and Vezzoni to the hyperbolic space in \cite{CV_indiana}, and a more unified approach to space forms $\mathbb M^n_+$  (the Euclidean space, the hyperbolic space and the hemisphere) has been done in \cite{CRV}, where more general functions of the principal curvatures have been considered. 

Let $S=\partial \Omega$ where $\Omega$ is a relatively compact connected open set in $M^n_+$,  and let $S$ be oriented by using the inward normal vector field to $\Omega$. Let $\{\kappa_1,\dots \kappa_{n-1}\}$ be the principal curvatures of $S$ ordered increasingly.  
Let ${\sf H}_S$ be one of the following functions: 
\begin{enumerate}
\item[i)] the mean curvature $H:=\tfrac{1}{n-1}\sum_{i}\kappa_i$;

\item[ii)] $f(\kappa_1,\dots,\kappa_{n-1})$, where 
$$
f\colon \{x=(x_1,\dots,x_{n-1})\in \mathbb \RR^{n-1}\,\,:\,\,x_1\leq x_2\leq\dots\leq x_{n-1} \}\to \RR\,,
$$ 
is a $C^2$-function such that 
$$
f(x)>0, \mbox{ if } x_i>0 \mbox{ for every }i=1,\dots,n-1
$$
and $f$ is concave on the component $\Gamma$ of $\{x\in\RR^{n-1}\,\,:\,\,f(x)>0\}$ containing  $\{x\in \RR^{n-1}\,\,:\,\, x_i>0\}$. 
\end{enumerate}   
For this class of functions, we can give quantitative estimates of proximity to a single ball in the spirit of \cite{CV_JEMS}.

\begin{theorem}\label{main_CRV}
Let $S$ be a $C^2$-regular, connected, closed hypersurface embedded in $\mathbb M^n_{+}$ satisfying a uniform touching ball condition of radius $\rho$. There exist constants $\ep,\, C>0$ such that if 
\begin{equation}\label{H quasi const}
{\osc}({\sf H}_S)
 \leq \ep,
\end{equation}
then there are two concentric balls $ B^d_{r}$ and $B^d_{R}$ of $\mathbb M^n_+$ such that
\begin{equation}\label{Bri Om Bre}
S \subset \overline{B}^{\,d}_{R} \setminus B^d_{r},
\end{equation}
and
\begin{equation}\label{stability radii}
R-r \leq C {\osc}({\sf H}_S).
\end{equation}
The constants $\ep$ and $C$ depend only on $n$ and upper bounds on $\rho^{-1}$ and on the area of $S$.

Moreover, $S$ is $C^{1}$-close to a sphere and there exists a $C^2$-regular map $\Psi: \partial B_r^d\to \mathbb{R}$ such that
$$
F(p) = \exp_x(\Psi(p)N_p) 
$$
defines a $C^2$-diffeomorphism from $B^d_r$ to $S$ and
\begin{equation}\label{fristLipbound}
\|\Psi\|_{C^1(\partial B_r)} \leq C\, {\osc}({\sf H}_S)^{1/2}\,,
\end{equation}
where $N$ is a normal vector field to $\partial B^d_r$. 
\end{theorem}

We notice that if $H_r$ denotes the $r$-higher order curvature of $S$ defined as the elementary symmetric polynomial of degree $r$ in the principal curvatures of $S$, then $H_r^{1/r}$ satisfies ii). 

Another important remark is to notice that the extension to non-Euclidean spaces requires to introduce a new concept of closedness between two hypersurfaces. In particular, in this more general context, the key estimate \eqref{bound on dist} is replaced by the following statement:

 for any $p \in\Sigma$ there exists $\hat p\in \hat\Sigma$ such that
\begin{equation*}
d(p,\hat p) +  |N_p-\tau_{\hat p}^p (N_{\hat p})|_p \leq C\, {\osc}({\sf H}_S)  \,,
\end{equation*}
where  $\tau_p^q\colon \RR^n\to \RR^n$ denotes the parallel transport along the unique geodesic path in $\mathbb H^n$ connecting $p$ to $q$.

\medskip

\subsection{Nonlocal Alexandrov theorem}  

In recent years, a lot of attention has been given to nonlocal problems. Starting from the definition of nonlocal perimeter \cite{CRS,CSi}
$$
P_s(\Om)=\int_\Om\int_{\Om^c}\frac{dx\,dy}{|x-y|^{n+2\,s}}\,,\qquad \Om^c = \mathbb R^n \setminus \Om\,,
$$
for $s\in(0,1/2)$, boundaries of sets $\Om\subset\mathbb R^n$ which are stationary for the $s$-perimeter functional are such that the nonlocal mean curvature is constant where, if $\partial \Omega$ is sufficiently smooth, the nonlocal mean curvature is given by
\begin{equation} \label{def:Hs}
H_s^\Omega (p) = \frac{1}{\omega_{n-2}}
\int_{\mathbb R^n} \frac{\widetilde\chi_\Om(x)}{|x-p|^{n+2s}} \,dx \,,\qquad \widetilde\chi_\Om(x)=\chi_{\Om^c}(x)-\chi_\Om(x)\,,
\end{equation}
where
$\chi_E$ denotes the characteristic function of a set $E$, $\omega_{n-2}$ is the measure of the $(n-2)$-dimensional sphere, and the integral is defined in the principal value sense. We notice that the nonlocal mean curvature can also be computed as a boundary integral, that is
 \begin{equation} \label{def:Hsdiv}
 H_s^\Om (p) = \frac{1}{s\, \omega_{n-2}} \int_{\pa \Om} \frac{(x-p)\cdot \nu_x}{|x-p|^{n+2s}}\, d \mathcal H^{n-1}_x \,.
 \end{equation}
where $\nu_x$ denotes the exterior unit normal to $\Om$ at $x\in \pa\Om$.
In \cite{Cabre} and \cite{CFMN}, the nonlocal version of the Alexandrov's 
theorem was proved:

\begin{theorem}\label{thm_nonlocal}
If $\Omega$ is a bounded open set of class $C^{1,2s}$ and $H_s^\Om$ is constant on $\pa\Om$, then $\pa\Om$ is a sphere.
\end{theorem}

The quantitative version of Theorem \ref{thm_nonlocal} was investigated in \cite{CFMN} where it is proved that if $H^\Omega_s$ has small Lipschitz constant then $\partial\Omega$ is close to a sphere, with a sharp estimate in terms of the deficit.
 
It is important to notice that the connectedness of $\Omega$ does not play a role here. This was indeed expected, since two disjoint balls does not have constant nonlocal mean curvature, since every point of $\partial \Omega$ influences the value of the nonlocal mean curvature at any other point of $\partial \Omega$. Hence, in the nonlocal case, one can not expect to have bubbling and the problem is more rigid.

As done in the local case, under suitable regularity assumption on $\partial \Omega$, in \cite{CFMN} we proved that $\partial \Omega$ is $C^{1,\alpha}$ close to a single sphere. Moreover, in the nonlocal case, we can prove something more, namely the $C^{2,\alpha}$ proximity to a single sphere (see \cite[Theorem 1.5]{CFMN}). This result gives an intriguing feature of the nonlocal case, which is the following: if the deficit is small then $\Omega$ is convex (and close to a single sphere).

\bibliographystyle{alpha}

\end{document}